\documentclass[10pt]{article}
\font\smallit=cmti10

\usepackage{amssymb,latexsym,amsmath,epsfig,amsthm} 
\makeatletter

\renewcommand\section{\@startsection {section}{1}{\z@}
{-30pt \@plus -1ex \@minus -.2ex}
{2.3ex \@plus.2ex}
{\normalfont\normalsize\bfseries\boldmath}}

\renewcommand\subsection{\@startsection{subsection}{2}{\z@}
{-3.25ex\@plus -1ex \@minus -.2ex}
{1.5ex \@plus .2ex}
{\normalfont\normalsize\bfseries\boldmath}}

\renewcommand{\@seccntformat}[1]{\csname the#1\endcsname. }

\makeatother

\newtheorem{theorem}{Theorem}
\newtheorem{lemma}{Lemma}
\newtheorem{proposition}{Proposition}

\theoremstyle{definition}
\newtheorem{example}{Example}
\newtheorem{definition}{Definition}

\begin{document}

\begin{center}
\uppercase{\bf The repetends of reduced fractions $\bf{\lowercase{a/b}^{\lowercase{k}}}$ approach\\full complexity with an increasing ${\lowercase{k}}$}
\vskip 20pt
{\bf Josefina L\'opez}\\
{\smallit Camataquí Villa Abecia,\\
Sud Cinti, Chuquisaca, Bolivia} \\
josefinapedro@hotmail.com
\vskip 10pt
{\bf Peter Stoll}\footnote{The majority of this work was done prior to the death of the Dr. phil. nat. Peter Stoll. The first author dedicates this paper with love, respect and admiration to his memory.}\\
{\smallit Camataquí Villa Abecia,\\
Sud Cinti, Chuquisaca, Bolivia}\\
\end{center}
\vskip 20pt

\vskip 30pt

\centerline{\bf Abstract}
\noindent

In this paper, we prove a criterion for complexity in $g$-ary expansions of a rational fraction $a/b<1$ with gcd$(a,b)=1$. We prove that for any purely periodic proper fraction $a/b$ and all $j\geq 1$, each sequence of $j$ digits occurs in the $g$-ary repetend of $a/b^k$ with a relative frequency that approaches $1/g^j$ with an increasing $k$. The absolute frequencies can be calculated by means of a simple transition matrix. Let $(a_k)$  be a sequence of positive integers relatively prime to $b$. We prove that each sequence of $j$ digits occurs in the $g$-ary repetend of $a_k/b^k$ with a relative frequency that approaches $1/g^j$ with an increasing $k$, unless all prime factors of $b$ divide the base $g\geq 2$.

\pagestyle{myheadings}
\thispagestyle{empty}
\baselineskip=12.875pt
\setcounter{page}{1}
\vskip 30pt


\section{Introduction}
A {\it proper fraction} is a rational number $a/b\in \mathbb{Q}$ such that $0<a/b<1$, where $a\in \mathbb{N}$, $b \in \mathbb{N}\setminus\{1\}$, and gcd$(a,b)=1$. Let $g\geq 2$ be a positive integer. Any proper fraction $a/b$ with a denominator $b$ relatively prime to $g$ has a unique $g${\it -ary expansion}
\begin{displaymath}
a/b=x_1g^{-1}+x_2g^{-2}+x_3g^{-3}+\cdots,
\end{displaymath}
where the digits $x_i$ for $i \in \mathbb{N}$ belong to the set $\{0,1,\ldots ,g-1\}$. We write this expansion
\begin{displaymath}
0.\,x_1x_2x_3\cdots .
\end{displaymath}

This $g$-ary expansion is purely periodic since gcd$(b,g)=1$. The smallest positive integer $n\in \mathbb{N}$ such that $x_i=x_{i+n}$ for all i $\in \mathbb{N}$ is the {\it period}, and the word $x_1x_2\cdots x_n$ is the {\it repetend} of the expansion.

Conversely, a {\it purely periodic} $g$-ary expansion $0.\,x_1x_2x_3\cdots$ represents a proper fraction $a/b$ with a denominator $b$ relatively prime to $g$, with the only exceptions $0=0.000\cdots$ and $1=0.(g-1)(g-1)(g-1)\cdots$. \\
The period of a purely periodic proper fraction $a/b$ is given by
\begin{displaymath}
\textrm{ord}_g(b):= \textrm{min} \{n\in \mathbb{N}\mid g^n\equiv 1 \pmod{b} \}.\footnote{Stoneham \cite{1969} defined $\omega (m):=$ ord$_g(b)$. }
\end{displaymath}
\begin{definition}\label{Fac}
\normalfont (Lothaire \cite{22002}). A word $f$ is called a {\it factor} of a word $x$ if there exist words $u$, $v$ such that $x=ufv$.
\end{definition}
\begin{definition}\label{F}
\normalfont Let $j\geq 1$ be a positive integer. The {\it enlarged repetend} of a purely periodic proper fraction $a/b$ is the finite word
\begin{displaymath}
F_g(\frac{a}{b},j):=(x_1x_2\cdots x_ {\textrm{ord}_g(b)})(x_1x_2\cdots x_{j-1}).
\end{displaymath}
\end{definition}
For $j=1$, the second factor is the empty word.

Since the $g$-ary expansion is purely periodic, any factor (subword) of length $j$ in the expansion is already a factor of $F_g(a/b,j)$.
\begin{definition}\label{m-word}
\normalfont A word of length $j\geq 1$ is called a {\it $j$-word}. The $j$-word at $i$ is the word
\begin{displaymath}
x_{i-j+1}x_{i-j+2}\cdots x_i,\quad\textrm{where }1\leq j\leq i.
\end{displaymath}
\end{definition}

The sequence of $j$-words at $i$ in the $g$-ary expansion of the purely periodic proper fraction $a/b$ is the sequence
\begin{displaymath}
\big(x_{i-j+1}x_{i-j+2}\cdots x_i\big)_{i=j}^{\infty}.
\end{displaymath}
However, any of its $j$-words can already be found in the finite sequence
\begin{displaymath}
S_g(\frac{a}{b},j):=\big(x_{i-j+1}x_{i-j+2}\cdots x_i\big)_{i=j}^{\textrm{ord}_g(b)+j-1}.
\end{displaymath}

Each $j$-word of $S_g(a/b,j)$ can be interpreted as a {\it base $g$ integer} with $j$ digits:
\begin{displaymath}
x_{i-j+1}x_{i-j+2}\cdots x_i = x_{i-j+1}g^{j-1}+x_{i-j+2}g^{j-2}+\cdots +x_ig^0,
\end{displaymath}
where $x_i$ is the least significant digit and where the leading digits are allowed to be zeros.

It is convenient to use base 10 numbers when referring to $j$-words of $S_g(a/b,j)$. Each $j$-word receives a {\it list number} $s$ from $0$ to $g^j-1$ which is its value when read as a base $g$ integer. The $j$-word at $i$ receives the list number
\begin{displaymath}
s=\sum_{m=0}^{j-1} x_{i-m}g^m.
\end{displaymath}
For instance, if $g=2$ and $j=3$, the list number of the word $110$ is $6$, and the list number of $011$ is $3$.

We denote the set of list numbers by
\begin{displaymath}
J_g:=\{0,1,\ldots ,g^j-1\}.
\end{displaymath}

Throughout this paper we  shall use the letter $s$ to represent either $j$-words or their respective list numbers. When translating a list number back into a word, we take care of eventually leading zeros. For instance, if $g=2$ and $j=4$, the list number $5$ is the word $0101$.

The {\it absolute frequency} $\nu(s)$\footnote{Stoneham \cite{1969} defined $N(B_j,g):=\nu(s)$, where $B_j:=s$. } is a function
\begin{displaymath}
\nu: J_g \longmapsto \mathbb{N}\cup  \{0\}
\end{displaymath}
that counts the number of times a $j$-word with list number $s\in J_g$ occurs as a factor in the enlarged repetend $F_g(a/b,j)$. There are exactly ord$_g(b)$ factors of length $j$ in $F_g(a/b,j)$ and therefore,
\begin{displaymath}
\sum_{s\in J_g} \nu(s) = \textrm{ord}_g(b).
\end{displaymath}

The {\it oscillation of absolute frequencies} is a nonnegative integer defined by
\begin{displaymath}
\sigma_g(a/b,j):=\textrm{max}\big\{|\nu(s)-\nu(s')|\,\big |\,s,s'\in J_g\big\}.
\end{displaymath}
The {\it relative frequency} of a $j$-word in $F_g(a/b,j)$ is defined by
\begin{displaymath}
\nu_{\textrm{rel}}(s):= \frac{\nu(s)}{\textrm{ord}_g(b)}.
\end{displaymath}
Of course,
\begin{displaymath}
\sum_{s\in J_g}\nu_{\textrm{rel}}(s)=1.
\end{displaymath}

\begin{definition}\label{Comp}
\normalfont Let $x_1x_2x_3\cdots$ be the binary expansion of $x$. The {\it complexity function} of $x$ counts the number $P(x,m)$ of different factors of length $m$ in $x$ for each nonnegative integer $m\geq0$.
\end{definition}

The number $P(x,1)$ counts the different letters appearing in $x$. In our case, the alphabet is $\{0,1\}$, so if $x$ is nonperiodic, $P(x,1) = 2$. On the other hand, if $x$ is purely periodic with a repetend of length $m$, we have $P(x,m) = m$.

In a former paper \cite{2012}, we studied the effect of successive divisions by $3$ on the complexity of the binary expansion of $1/3^k$. We generalize these results for all reduced fractions of the form $a/b^k$ in any integer base $g\geq 2$.

Following Waldschmidt \cite{2009}, a real number $\xi$ is called {\it simply normal in base} $g$ if each digit occurs with (well-defined limiting relative) frequency $1/g$ in its $g$-ary expansion. Furthermore, a real number $\xi$ is called {\it normal in base} $g$ if it is simply normal in base $g^j$ for all $j\geq 1$. Hence a real number $\xi$ is normal in base $g$ if and only if, for all $j\geq 1$, each sequence of $j$ digits occurs with (limiting relative) frequency $1/g^j$ in its $g$-ary expansion.

Our issue is far simpler. The fractions have an eventually periodic $g$-ary expansion. Nevertheless, we borrow the concept of normality for our purpose since it describes exactly the digits structure of an increasingly long repetend under successive divisions by $b$. With an increasing $k$, larger than an initial constant $c_g(b)$, which depends only on $b$ and $g$, the $g$-ary period ord$_g(b^k)$ is strictly increasing unless all prime factors of $b$ divide $g$.

We prove that for any purely periodic proper fraction $a/b$ and all $j\geq 1$, each sequence of $j$ digits occurs in the $g$-ary repetend of $a/b^k$ with a relative frequency that approaches $1/g^j$ with an increasing $k$ (Theorem \ref{pure}). This proof is given in Sections \ref{SeqW} and \ref{Matrix}.

Also, in Section  \ref{SeqW}, we prove that the set of (distinct) $j$-words in $F_g(a/b^{k+t},j)$ is completely determined by the corresponding set of $j$-words in $F_g(a/b^k,j)$ for all sufficiently large $k$ and $t\in \mathbb{N}$ (Proposition \ref{Xkt}). Next, in Section  \ref{Crit}, we prove a criterion for complexity in g-ary expansions of rational fractions (Theorem \ref{Criterion}). Further, in Section  \ref{Matrix}, we show that the absolute frequencies can be calculated by means of a simple transition matrix (Proposition \ref{Dt} and Example \ref{Ex}), which guarantees that the oscillation of absolute frequencies of $j$-words in $F_g(a/b^k,j)$ is bounded for all $k\in \mathbb{N}$ (Proposition \ref{cota} and Example \ref{Ex2}). Finally, in Section \ref{Prop}, we generalize Theorem \ref{pure} for all reduced fractions of the form $a_k/b^k$ (Proposition \ref{seq} and Theorem \ref{proper}). 


\section{Sequences of Words and Remainders}\label{SeqW}
\begin{lemma}\label{qi}
Let $0.\,x_1x_2\cdots$ be the $g$-ary expansion of a purely periodic proper fraction $a/b$. For all $i\in \mathbb{N}$, the quotient $q_i= \Big \lfloor ag^i/b \Big \rfloor$
is the $g$-ary integer $q_i=x_1x_2\cdots x_i$.
\end{lemma}
\noindent {\it Proof.} 
(a) We have $a/b=x_1g^{-1}+x_2g^{-2}+\cdots$. Thus
\begin{eqnarray*}
\frac{ag^i}{b}&=&x_1g^{i-1}+ x_2g^{i-2}+x_3g^{i-3}+\cdots+ x_ig^{i-i}+x_{i+1}3g^{i-{(i+1)}}+\cdots\\
&<& x_1g^{i-1}+ x_2g^{i-2}+\cdots+ x_ig^{i-i}+(g-1)g^{-1}+(g-1)g^{-2}+\cdots\\
&=&x_1g^{i-1}+ x_2g^{i-2}+\cdots+ x_ig^0+1
\end{eqnarray*}
imply that
\begin{displaymath}
q_i=\Big\lfloor\frac{ag^i}{b}\Big\rfloor.
\end{displaymath}
Hence $q_i=x_1x_2\cdots x_i$.\hfill $\Box$
\begin{lemma} \label{Sg}
The sequence of $j$-words at $i$ in the enlarged $g$-ary repetend $F_g(a/b,j)$ of a purely periodic proper fraction is given by the following sequence of list numbers:
\begin{displaymath}
S_g(\frac{a}{b},j)=\Big ( \Big\lfloor\frac{ag^i}{b}\Big\rfloor \bmod{g^j}\Big )_{i=j}^{\mathrm{ord}_g(b)+j-1}.
\end{displaymath}
\end{lemma}
\noindent {\it Proof.} By Lemma \ref{qi}, the integer quotient $q_i=\lfloor ag^i/b\rfloor$ is the $g$-ary integer $x_1x_2\cdots x_i$. If $i\geq j$, then
\begin{displaymath}
(x_1x_2\cdots x_i) \bmod{g^j}=x_{i-j+1}x_{i-j+2}\cdots x_i.\\[-18pt]
\end{displaymath}\hfill $\Box$

\noindent Note that $S_g(a/b,1)$ is the repetend of $a/b$.

For a fixed $b\not=1$ relatively prime to $g$, the set of purely periodic proper fractions with denominator $b$ splits in exactly $\varphi(b)/\textrm{ord}_g(b)$ equivalence classes, where $\varphi$ is Euler's totient function: any two of these fractions are $g$-equivalent if their numerators belong to the same coset of the cyclic subgroup $D\subset\big(\mathbb{Z}/b\mathbb{Z}\big )^\times $ generated by $g$ modulo $b$. Indeed, each numerator $a$ relatively prime to $b$ defines the coset
\begin{displaymath}
aD=\big\{ag^i\,\big | \,a\in\big(\mathbb{Z}/b\mathbb{Z}\big )^\times ,\,\,\,  g^i\in D \quad\textrm{for all } i\in\{1,2,\ldots , \textrm{ord}_g(b)\}\big\}.
\end{displaymath}
We define the {\it rotate left operator} $\rho$ by
\begin{displaymath}
\rho (s_1s_2s_3\cdots s_{t}):=s_2s_3\cdots s_{t}s_1.
\end{displaymath}
The following lemma shows that $g$-equivalent purely periodic proper fractions have up to successive left rotations the same $j$-words in the enlarged repetend.
\begin{lemma}\label{rho}
Let $a/b$ be a purely periodic proper fraction and $j,n\in\mathbb{N}$. Then we have
\begin{displaymath}
S_g\Big (\frac{(ag^n)\bmod{b}}{b},j\Big )=\rho ^n\Big (S_g(\frac{a}{b},j)\Big ).
\end{displaymath}
\end{lemma}
\noindent {\it Proof.} For $n=1$,
we have $a/b=x_1x_2\cdots$. Hence $ag/b=x_1.\,x_2x_3\cdots$. Therefore,
\begin{displaymath}
\frac{(ag)\bmod{b}}{b}=0.\,x_2x_3\cdots .
\end{displaymath}
By comparing the two sequences
\begin{eqnarray*}
S_g(a/b,j)&=&(x_1x_2\cdots x_j,\,x_2x_3\cdots x_{j+1},\ldots ,\,x_{\textrm{ord}_g(b)}x_1x_2\cdots x_{j-1}),\\
S_g\Big(\big ((ag)\bmod{b}\big )/b,j\Big)&=&(x_2x_3\cdots x_{j+1},\ldots ,\,x_{\textrm{ord}_g(b)}x_1x_2\cdots x_{j-1},\,x_1x_2\cdots x_j),
\end{eqnarray*}
we get the statement for $n=1$. For any other $n\in\{2,3,\ldots ,\textrm{ord}_g(b)\}$ a similar argument holds by comparing the following expansions:
\begin{eqnarray*}
\big((ag^{n-1})\bmod{b}\big)/b&=&0.\,x_{n}x_{n+1}\cdots ,\\
\big((ag^n)\bmod{b}\big)/b&=&0.\,x_{n+1}x_{n+2}\cdots .\\[-36pt]
\end{eqnarray*}\hfill $\Box$

\begin{lemma}\label{rs}
Let $a/b$ be a purely periodic proper fraction, $i\in \mathbb{N}$, $j\in\mathbb{N}$ and $r_i:=(ag^i) \bmod{b}$. For all $i\geq j$, the $(i-j)$-th remainder
\begin{displaymath}
r_{i-j}:= (ag^{i-j})\bmod{b}
\end{displaymath} 
of the long division $a$ by $b$ determines the list number $s_i$ of the $j$-word at $i$ by the equation
\begin{displaymath}
r_{i-j} \frac{g^j}{b}=s_i+\frac{r_i}{b},
\end{displaymath}
and consequently,
\begin{displaymath}
\Big\lfloor r_{i-j}\frac{g^j}{b}\Big\rfloor = s_i.
\end{displaymath}
\end{lemma}
\noindent {\it Proof.} For $i\geq j$, let
\begin{eqnarray}
ag^i &=& q_ib+r_i \quad (0<r_i<b),\label{eins}\\ 
ag^{i-j} &=& q_{i-j}b+r_{i-j}\quad (0<r_{i-j}<b),\label{zwei}\\
q_i &=& c_ig^j+s_i \quad (0<s_i<g^j).\label{drei}
\end{eqnarray}
Note that $r_0\equiv a \bmod{b} = a$ by Equation (\ref{zwei}), hence $r_0 = a$.

By Lemma \ref{Sg} and Equation (\ref{eins}), $s_i$ is the list number of the $j$-word at $i$. We multiply Equation (\ref{zwei}) by $g^j$ and substitute this expression together with Equation (\ref{drei}) into Equation (\ref{eins}). We get
\begin{equation}
g^j(q_{i-j}b+r_{i-j})=(c_ig^j+s_i)b+r_i. \label{vier}
\end{equation}

Furthermore, for $i>j$, we have the $g$-ary integers
\begin{eqnarray*}
q_{i-j}g^j&=&x_1x_2\cdots x_{i-j}00\cdots 0 \quad (j \textrm{ zeros}),\\
s_i&=&x_{i-j+1}x_{i-j+2}\cdots x_i.
\end{eqnarray*}
Hence $q_{i-j}g^j+s_i=x_1x_2\cdots x_i=q_i$. Comparing this expression with Equation (\ref{drei}), we get
\begin{equation} 
c_i=q_{i-j}.\label{fünf}
\end{equation}

This is even true for $i=j$, having now the $g$-ary integers
\begin{eqnarray*}
q_{i-j}g^j&=&00\cdots 0 \quad (j \textrm{ zeros}),\\
s_i&=&x_1x_2\cdots x_j,
\end{eqnarray*}
since $q_{i-j}=q_0=0$ by Equation (\ref{zwei}).

Finally, we substitute Equation (\ref{fünf}) into Equation (\ref{vier}). After dividing by $b$, we obtain the claimed relation.\hfill $\Box$

\begin{lemma}\label{ld}
Let $a/b$ be a purely periodic proper fraction and $j\in\mathbb{N}$. The remainders of the sequence
\begin{displaymath}
R_g(\frac{a}{b},j):=\Big ( (ag^{i-j})\bmod{b}\Big )_{i=j}^{\mathrm{ord}_g(b)+j-1}
\end{displaymath}
are all different from each other.
\end{lemma}
\noindent {\it Proof.} Let
\begin{displaymath}
R=\Big ( (ag^i)\bmod{b}\Big )_{i=0}^{\textrm{ord}_g(b)-1}.
\end{displaymath}
Of course, $R_g(a/b,j)$ and $R$ are the same sequence. The sequence $R$ gives the successive remainders of the long division $a$ by $b$, and ord$_g(b)$ is the period of the $g$-ary expansion.\hfill $\Box$

\begin{lemma}\label{h}
Let $b>1$ be an integer relatively prime to $g\geq 2$ and $k\in\mathbb{N}$. Define
\begin{displaymath}
h:=\frac{g^{\mathrm{ord}_g(b^k)}-1}{b^k}.
\end{displaymath}
If \, $d=\mathrm{gcd}(h,b)$,\, then \,$\mathrm{ord}_g(b^{k+1})=\frac{b}{d}\mathrm{ord}_g(b^k)$.
\end{lemma}
\noindent {\it Proof.} The integer $h$ is well-defined since
\begin{displaymath}
g^{\mathrm{ord}_g(b^k)}\equiv 1\pmod{b^k}.
\end{displaymath}
We consider two cases.

(a) If $b=d$, then
\begin{displaymath}
g^{\mathrm{ord}_g(b^k)}-1=hb^k=nb^{k+1} \quad\textrm{for some }n\in\mathbb{N}.
\end{displaymath}
Thus
\begin{displaymath}
g^{\mathrm{ord}_g(b^k)}\equiv 1\pmod{b^{k+1}},
\end{displaymath}
implying that ord$_g(b^k)$ is a multiple of ord$_g(b^{k+1})$. However, ord$_g(b^{k+1})$ is a multiple of ord$_g(b^k)$ since
\begin{displaymath}
g^{\mathrm{ord}_g(b^{k+1})}\equiv 1\pmod{b^{k+1}}\quad \textrm{implies that}\quad g^{\mathrm{ord}_g(b^{k+1})}\equiv 1\pmod{b^{k}}.
\end{displaymath}
Hence ord$_g(b^{k+1})=\,\,$ord$_g(b^k)$.

(b) If
\begin{displaymath}
g^{n\cdot\mathrm{ord}_g(b^k)}\equiv 1\pmod{b^{k+1}}
\end{displaymath}
for some integer $n>1$, then we get\footnote{The symbol ${\bf O}\big (hb^{k+1}\big )$ means that all the remaining terms are divisible by $hb^{k+1}$.}
\begin{displaymath}
g^{n\cdot\mathrm{ord}_g(b^k)}=1+n(hb^k)+{\bf O}\big (hb^{k+1}\big )
\end{displaymath}
by raising both sides of
\begin{displaymath}
g^{\mathrm{ord}_g(b^k)}=1+hb^k
\end{displaymath}
to the $n$-th power. It follows that $b^{k+1}$ must divide $nhb^k$. So we have
\begin{displaymath}
b\,\mid \, nh \quad \textrm{and}\quad\textrm{gcd}(h,b)=d.
\end{displaymath}
Hence $b/d$ divides $n$. The smallest of such $n$ is $n=b/d$ and furthermore,
\begin{displaymath}
g^{\frac{b}{d}\mathrm{ord}_g(b^k)}\equiv 1\pmod{b^{k+1}}
\end{displaymath}
since
\begin{displaymath}
\frac{b}{d}hb^k=\frac{h}{d}b^{k+1}
\end{displaymath}
is an integer divisible by $b^{k+1}$.\hfill $\Box$

The following lemma about the period of $a/b^k$ in base $g\geq 2$ is well-known. A sketched proof for $g=10$ can be found in Wikipedia.\footnote{Wikipedia, available at {\em http://en.wikipedia.org/wiki/Repeating\underline{ }decimal}.}

\begin{lemma}\label{bk} 
Let $b>1$ be an integer relatively prime to $g\geq 2$. There is some $c\in \mathbb{N}$ such that
\begin{displaymath}
\mathrm{ord}_g(b^{c+t})=b^t\mathrm{ord}_g(b^c)
\end{displaymath}
for all positive integers $t$.
\end{lemma}
\noindent {\it Proof.} The assertion follows from a lemma of Korobov \cite{1972}, available as Lemma 4.3 in Bailey and Crandall \cite{12002}.\hfill $\Box$
\begin{definition}\label{cg}
\normalfont Let $b>1$ be an integer relatively prime to $g\geq 2$. The integer $c_g(b)$ denotes the smallest $c\in\mathbb{N}$ satisfying Lemma \ref{bk}; that is,
\begin{displaymath}
c_g(b):= \mathrm{min}\, \{c\in \mathbb{N}\,\mid \, \mathrm{ord}_g(b^{c+t})=b^t\mathrm{ord}_g(b^c)\quad \textrm{for all }t\in\mathbb{N}\}.
\end{displaymath}
\end{definition}

\begin{table}[!htb]
\begin{center}
\begin{tabular}{|c|ccccc|c|ccccc|c|ccccc|}
\hline
&&&$g$&&&&&&$g$&&&&&&$g$&&\\
$b$&2&3&5&7&10&$b$&2&3&5&7&10&$b$&2&3&5&7&10\\
\hline
$2$&&3&2&4&&$35$&1&1&&&&$68$&&3&3&4&\\
$3$&1&&1&1&2&$36$&&&2&2&&$69$&1&&1&1&2\\
$4$&&2&1&2&&$37$&1&1&1&1&1&$70$&&4&&&\\
$5$&1&1&&2&&$38$&&3&2&4&&$71$&1&1&1&1&1\\
$6$&&&3&4&&$39$&2&&1&2&3&$72$&&&1&2&\\
$7$&1&1&1&&1&$40$&&2&&2&&$73$&1&1&1&1&1\\
$8$&&1&1&2&&$41$&1&1&1&1&1&$74$&&3&4&4&\\
$9$&1&&1&1&1&$42$&&&3&&&$75$&1&&&1&\\
$10$&&4&&5&&$43$&1&1&1&1&1&$76$&&2&1&2&\\
$11$&1&2&1&1&1&$44$&&2&1&2&&$77$&1&2&1&&1\\
$12$&&&2&2&&$45$&1&&&2&&$78$&&&4&5&\\
$13$&1&1&1&1&1&$46$&&3&3&4&&$79$&1&1&1&1&1\\
$14$&&3&3&&&$47$&1&1&1&1&1&$80$&&1&&2&\\
$15$&1&&&2&&$48$&&&1&1&&$81$&1&&1&1&1\\
$16$&&1&1&1&&$49$&1&1&1&&1&$82$&&5&4&6&\\
$17$&1&1&1&1&1&$50$&&4&&5&&$83$&1&1&1&1&1\\
$18$&&&3&4&&$51$&1&&1&1&2&$84$&&&2&&\\
$19$&1&1&1&1&1&$52$&&2&2&3&&$85$&1&1&&2&\\
$20$&&2&&3&&$53$&1&1&1&1&1&$86$&&3&3&4&\\
$21$&2&&2&&3&$54$&&&3&4&&$87$&1&&1&1&2\\
$22$&&3&2&4&&$55$&2&2&&3&&$88$&&2&1&2&\\
$23$&1&1&1&1&1&$56$&&1&1&&&$89$&1&1&1&1&1\\
$24$&&&1&2&&$57$&3&&3&2&4&$90$&&&&5&\\
$25$&1&1&&1&&$58$&&4&3&4&&$91$&1&1&1&&1\\
$26$&&3&4&5&&$59$&1&1&1&1&1&$92$&&2&2&2&\\
$27$&1&&1&1&1&$60$&&&&3&&$93$&1&&2&2&3\\
$28$&&2&2&&&$61$&1&1&1&1&1&$94$&&3&3&4&\\
$29$&1&1&1&1&1&$62$&&3&2&4&&$95$&1&1&&2&\\
$30$&&&&5&&$63$&1&&1&&2&$96$&&&1&1&\\
$31$&1&1&1&1&1&$64$&&1&1&1&&$97$&1&1&1&1&1\\
$32$&&1&1&1&&$65$&1&1&&2&&$98$&&3&3&&\\
$33$&1&&1&1&2&$66$&&&3&4&&$99$&1&&1&1&1\\
$34$&&6&6&7&&$67$&1&1&1&1&1&$100$&&2&&3&\\
\hline
\end{tabular}
\caption{Function $c_g(b)$ (Definition \ref{cg}) for $g\in\{2,3,5,7,10\}$, $2\leq b \leq 100$ and gcd$(g,b)=1$.}
\label{pgm}
\end{center}
\end{table}

\begin{lemma}\label{k1}
Let $a/b^k$ be a purely periodic proper fraction and $k\geq c_g(b)$. The remainders $r'_v$ in the sequence
\begin{displaymath}
R_g\Big (\frac{a}{b^{k+1}},j\Big )=\Big ( (ag^v) \bmod{b^{k+1}}\Big )_{v=0}^{\mathrm{ord}_g(b^{k+1})-1}
\end{displaymath}
have the form
\begin{displaymath}
r'_v=r_{i+m\cdot \mathrm{ord}_g(b^k)}=r_i+nb^k \quad (m,n\in \{0,1,\ldots ,b-1\}),
\end{displaymath}
where $r_i\in R_g(a/b^k,j)$. This representation is unique.
\end{lemma}
\noindent {\it Proof.} (a) By Lemma \ref{bk},
\begin{equation}\label{1}
h=\frac{g^{\textrm{ord}_g(b^k)}-1}{b^k}
\end{equation}
is relatively prime to $b$ since $k\geq c_g(b)$ and ord$_g(b^{k+1})=b$ ord$_g(b^k)$.\\
(b) For $0\leq i<$ ord$_g(b^k)$, we have
\begin{eqnarray}
ag^i&=&q_ib^k+r_i\quad (0<r_i<b^k),\nonumber\\
ag^i&=&q'_ib^{k+1}+r'_i\quad (0<r'_i<b^{k+1}).\label{2}
\end{eqnarray}

We write $ag^i$ as a $b$-adic sum; that is,
\begin{displaymath}
ag^i=y_0b^0+y_1b^1+\cdots +y_kb^k+\cdots ,
\end{displaymath}
where the $y$'s belong to the set $\{0,1,\ldots ,b-1\}$. Of course, only finitely many terms in this sum are not equal to zero.\footnote{As opposed to the g-ary expansions, the $b$-adic expansion of a real number is unique.} So we have
\begin{eqnarray*}
r_i&\equiv&(ag^i)\bmod{b^k}\quad =y_0b^0+\cdots +y_{k-1}b^{k-1},\\
r'_i&\equiv&(ag^i)\bmod{b^{k+1}}=y_0b^0+\cdots +y_{k-1}b^{k-1}+y_kb^k.
\end{eqnarray*}
Hence
\begin{equation}\label{3}
r'_i=r_i+y_kb^k.
\end{equation}
(c) Applying Equations (\ref{1}), (\ref{2}) and (\ref{3}), we get
\begin{eqnarray*}
ag^{i+\mathrm{ord}_g(b^k)}&=&(q'_ib^{k+1}+r_i+y_kb^k)(hb^k+1)\\
&=&r_i+(r_ih+y_k)b^k+{\bf O}\big (b^{k+1}\big ).
\end{eqnarray*}
We substitute $r_ih+y_k$ with the $b$-adic sum
\begin{displaymath}
r_ih+y_k:=z_0b^0+z_1b^1+\cdots\quad (z_i\in \{0,1,\ldots ,b-1\}).
\end{displaymath}
Since $z_0\equiv(r_ih+y_k)\bmod{b}$ we get
\begin{displaymath}
ag^{i+\textrm{ord}_g(b^k)}=r_i+\big ((r_ih+y_k)\bmod{b}\big )b^k+{\bf O}\big (b^{k+1}\big ).
\end{displaymath}
Hence
\begin{equation}\label{4}
r'_{i+\textrm{ord}_g(b^k)}=r_i+\big ((r_ih+y_k)\bmod{b}\big )b^k.
\end{equation}
(d) The method used in (b) and (c) leads to the following generalization of Equations (\ref{3}) and (\ref{4}):
\begin{equation}\label{5}
r'_{i+m\cdot\textrm{ord}_g(b^k)}=r_i+\big ((mr_ih+y_k)\bmod{b}\big )b^k
\end{equation}
for all $m\in\{0,1,\ldots ,b-1\}$.

The integers $r_i$ and $h$ are relatively prime to $b$ and therefore,
\begin{displaymath}
\Big ( (mr_ih+y_k)\bmod{b}\Big )_{m=0}^{b-1}
\end{displaymath}
is a permutation of $(0,1,\cdots,b-1)$. It follows that there are $b\cdot$ord$_g(b^k)$ integers on both sides of (\ref{5}). The remainders on the left are all different from each other by Lemma \ref{ld}.\hfill $\Box$
\begin{lemma}\label{kt}
Let $a/b^k$ be a purely periodic proper fraction, $t\in\mathbb{N}$ and $k\geq c_g(b)$. The remainders $r'_v$ in the sequence
\begin{displaymath}
R_g\Big (\frac{a}{b^{k+t}},j\Big )=\Big ( (ag^v) \bmod{b^{k+t}}\Big )_{v=0}^{\mathrm{ord}_g(b^{k+t})-1}
\end{displaymath}
have the form
\begin{displaymath}
r'_v=r_{i+m\cdot \mathrm{ord}_g(b^k)}=r_i+nb^k \quad (m,n\in \{0,1,\ldots ,b^t-1\}),
\end{displaymath}
where $r_i\in R_g(a/b^k,j)$. This representation is unique.
\end{lemma}
\noindent {\it Proof.} Since $k\geq c_g(b)$ we have $\textrm{ord}_g(b^{k+t})=b^t\textrm{ord}_g(b^k)$. We apply Lemma \ref{k1} successively. The exponent in the brackets of $i^{(k)}$ denotes the level at which the index $i$ is taken. For instance, the statement of Lemma \ref{k1} now looks like this:
\begin{displaymath}
v:=i^{(k+1)}=i^{(k)}+m_0\cdot \textrm{ord}_g(b^k)\quad (0\leq i^{(k)}<\textrm{ord}_g(b^k);\,0\leq m_0<b)
\end{displaymath}
and
\begin{displaymath}
r_{i^{(k)}+m_0\cdot \textrm{ord}_g(b^k)}=r_{i^{(k)}}+n_0b^k\quad (m_0,n_0\in\{0,1,\ldots ,b-1\}),
\end{displaymath}
where $r_{i^{(k)}}\in R_g(a/b^k,j)$.

For the index $i^{(k+t)}$, we get the following:
\begin{eqnarray*}
i^{(k+2)}&=&i^{(k+1)}+m_1\cdot \textrm{ord}_g(b^{k+1})\quad (0\leq i^{(k+1)}<b\cdot\textrm{ord}_g(b^k);\,0\leq m_1<b)\\
&=&i^{(k)}+(m_0+m_1b)\cdot \textrm{ord}_g(b^k),\\
i^{(k+3)}&=&i^{(k+2)}+m_2\cdot \textrm{ord}_g(b^{k+2})\quad (0\leq i^{(k+2)}<b^2\cdot\textrm{ord}_g(b^k);\,0\leq m_2<b)\\
&=&i^{(k)}+(m_0+m_1b+m_2b^2)\cdot \textrm{ord}_g(b^k),\\
&\vdots&
\end{eqnarray*}
\begin{equation*}
i^{(k+t)} = i^{(k)}+\Big(  \sum_{x=0}^{t-1}m_xb^x  \Big) \textrm{ord}_g(b^k)\quad (0\leq i^{(k+t-1)}<b^{t-1}\cdot\textrm{ord}_g(b^k);\,0\leq m_x<b).
\end{equation*}
Let $m:=\sum_{x=0}^{t-1}m_xb^x$. Then we have $0\leq m\leq b^t-1$ since this $b$-adic sum yields exactly $t$ different nonnegative integers less than $b^t$.

For the remainder $r_{i^{(k+t)}}$, we get the following:
\begin{eqnarray*}
r_{i^{(k+1)}}&=&r_{i^{(k)}}+n_0b^k\quad (0\leq n_0<b),\\
r_{i^{(k+2)}}&=&r_{i^{(k+1)}}+n_1b^{k+1}\quad (0\leq n_1<b)\\
&=&r_{i^{(k)}}+(n_0+n_1b)b^k,\\
&\vdots&\\
r_{i^{(k+t)}}&=&r_{i^{(k)}}+\Big (\sum_{x=0}^{t-1}n_xb^x\Big )b^k\quad (0\leq n_x<b).
\end{eqnarray*}
Hence $n:=\sum_{x=0}^{t-1}n_xb^x$ is the claimed $n$.\hfill $\Box$
\begin{definition}\label{de}
\normalfont Let $X_g(a/b^k,j)$ denote the set of (distinct) list numbers of $j$-words occurring in the enlarged $g$-ary repetend $F_g(a/b^k,j)$.
\end{definition}
\begin{proposition}\label{Xkt}
Let $a/b^k$ be a purely periodic proper fraction, $k\geq c_g(b)$ and $j,t\in\mathbb{N}$. We have
\begin{displaymath}
X_g(\frac{a}{b^{k+t}},j)=\Big\{ \Big ( \frac{1}{b^t}(s-n)\Big )\bmod{g^j}\, \Big | \, s\in X_g(\frac{a}{b^k},j) \,\, \textrm{and } n\in\big\{0,1,\ldots, b^t-1\big\} \Big\},
\end{displaymath}
where $c_g(b)$ is the integer of Definition \ref{cg}.
\end{proposition}
\noindent {\it Proof.} Let
\begin{eqnarray*}
R_g(\frac{a}{b^{k+t}},j)&=&\Big ((ag^{v-j})\bmod{b^{k+t}}\Big )_{v=j}^{\textrm{ord}_g(b^{k+t})+j-1},\\
R_g(\frac{a}{b^k},j)&=&\Big ((ag^{i-j})\bmod{b^k}\Big )_{i=j}^{\textrm{ord}_g(b^k)+j-1}.
\end{eqnarray*}
Of course, the first sequence is the same sequence as in Lemma \ref{kt}, and the second one is given in Lemma \ref{ld}.

Each $r'_v\in R_g(a/b^{k+t},j)$ has the form
\begin{equation}\label{uno}
r'_v=r_{i-j}+nb^k,
\end{equation}
where $v=i-j+m\cdot\textrm{ord}_g(b^k)$, $j\leq i\leq \textrm{ord}_g(b^k)+j-1$ and $m,n\in \{0,1,\ldots ,b^t-1\}$.

By Lemma \ref{rs}, we have
\begin{displaymath}
r'_{v-j}\frac{g^j}{b^{k+t}}=s'_v+\frac{r'_v}{b^{k+t}}\quad\quad\quad \big(s'_v\in  S_g(a/b^{k+t},j)\big).
\end{displaymath}
Hence
\begin{equation}\label{dos}
-\frac{1}{b^{k+t}}r'_v\equiv s'_v\pmod{g^j}.
\end{equation}
Equations (\ref{uno}) and (\ref{dos}) yield
\begin{displaymath}
-\frac{1}{b^{k+t}}(r_{i-j}+nb^k)\equiv s'_v\pmod{g^j}.
\end{displaymath}
But by Lemma \ref{rs},
\begin{displaymath}
-\frac{1}{b^k}r_{i-j}\equiv s_{i-j}\pmod{g^j}\quad \quad\quad\big(s_{i-j}\in  S_g(a/b^k,j)\big)
\end{displaymath}
and therefore,
\begin{equation}\label{tres}
\frac{1}{b^t}(s_{i-j}-n)\equiv s'_v\pmod{g^j}.
\end{equation}

For any $s'_v\in S_g(a/b^{k+t},j)$ --- and hence for any $j$-word of $X_g(a/b^{k+t},j)$ --- there is a $s_{i-j}\in S_g(a/b^k,j)$ such that Equation (\ref{tres}) holds. Indeed, $v=i-j+m\cdot\textrm{ord}_g(b^k)$ implies that $i-j$ is equal to $v$ modulo ord$_g(b^k)$.\hfill $\Box$


\section{A Criterion for Complexity}\label{Crit}
\begin{definition}\label{S}
\normalfont
Let $a/b<1$ be a rational fraction in the lowest terms and $j\geq1$. The {\it set of $j$-words} occurring in the expansion of $a/b$ is the finite set
\begin{displaymath}
S_R(\frac{a}{b},j):=\{\lfloor r_ig^j/b\rfloor |i=0,1,2,\ldots \},
\end{displaymath}
where $r_i$ is a power residue.
\end{definition}

Let $a/b$ be a rational fraction in lowest terms. If
\begin{itemize}
\item $b=b'\prod_{p_i|g}p_i^{b_i}$ is such that gcd$(g,b')=1$,
\item $N$ is the number of digits in the non-periodic part,
\end{itemize}
we define
\begin{displaymath}
dr :=\left\{\begin{array}{ll} 
N+1&\textrm{if }b'=1;\\
N+\textrm{ord}_g(b')&\textrm{else}.
\end{array}\right .
\end{displaymath}
This definition guarantees that for a fixed $j\geq 1$, the words $s\in S_R(\frac{a}{b},j)$ have their first digit $x_{i-j+1}$ somewhere in
\begin{displaymath}
.x_1x_2\cdots x_Nx_{N+1}\cdots x_{dr}
\end{displaymath}
and any other $j$-word in the $g$-ary expansion of $a/b$ is equal to some $s\in S_R(\frac{a}{b},j)$.
\begin{definition}\label{ts}
\normalfont
Let $g\geq 2$ and $j\geq 1$ be fixed integers such that $0<g^j/b\leq 1$, and let
\begin{displaymath}
J_g=\{0,1,\ldots ,g^j-1\},\quad B=\{0,1,\dots ,b-1\}.\footnote{Note that $b\geq 2$.}
\end{displaymath}
For every fixed $s\in J_g$, we define
\begin{displaymath}
t_s:=\mathrm{min}\Big\{t\in B\Big |\Big\lfloor t\cdot \frac{g^j}{b}\Big\rfloor = s\Big\}
\end{displaymath}
and say that
\begin{displaymath}
T:=\Big (\big (t_s\big )_{s=0}^{g^j-1},b\Big ):=\big (t_s\big )_{s=0}^{g^j}
\end{displaymath}
is the {\it master sequence}.
\end{definition}
By Lemma \ref{generalizar}, the master sequence is well-defined. It is a strictly increasing sequence which gives the position $t_s\in (0,1,\ldots ,b-1,b)$ of the first $s$ in a block of $\lceil b/g^j\rceil$ or $\lfloor b/g^j\rfloor$ repeating integers $s$ in $(\lfloor tg^j/b \rfloor)_{t=0}^{b-1}$. For all $s\in J_g$, we have
\begin{displaymath}
\lfloor b/g^j\rfloor\leq t_{s+1}-t_s\leq \lceil b/g^j\rceil.
\end{displaymath}

From the inequality above, and since $s=\lfloor r_{i-j}g^j/b \rfloor$ if and only if $t_s\leq r_{i-j}<t_{s+1}$, with $i\geq j$, it follows that any given $j$-word $s=x_{i-j+1}\cdots x_i$ occurs in the first $d+j-1$ digits of the $g$-ary expansion not more than $\lceil b/g^j\rceil$ times.
The master sequence depends on $g$, $b$ and $j$. Luckily, its terms can be explicitly calculated.
\begin{lemma}\label{T}
Let $g\geq 2$ and $j\geq 1$ be fixed integers such that $0<g^j/b\leq 1$. The master sequence is given by
\begin{displaymath}
T=\bigg(\frac{sb+\big ( (sh_1)\bmod{g^j}\big )}{g^j}\bigg )_{s=0}^{g^j},
\end{displaymath}
\begin{displaymath}
\textrm{where}\quad h_1:=\Big (\Big\lceil\frac{b}{g^j}\Big\rceil g^j\Big )\bmod{b}.
\end{displaymath}
\end{lemma}
\noindent {\it Proof.} Of course, $t_0=0$ and $t_{g^j}=b$. By Lemma \ref{generalizar}, the number $0$ occurs $\lceil b/g^j\rceil$ times. Therefore, $t_1=\lceil b/g^j\rceil$. Let us write
\begin{displaymath}
t_sg^j=\lfloor t_sg^j/b\rfloor b+h_s\quad (0\leq h_s<b).
\end{displaymath}
Thus
\begin{displaymath}
t_sg^j=sb+h_s\quad\textrm{and}\quad t_1g^j=b+h_1.
\end{displaymath}
But $st_1g^j=sb+sh_1$ yields $h_s \equiv (sh_1)\bmod{g^j}$, and hence $t_sg^j=sb+h_s=sb+(sh_1)\bmod{g^j}$.\hfill $\Box$

Let us say that a $g$-ary expansion has {\it complexity} $g^j$ if all the $g^j$-words of fixed length $j\geq1$ are factors (sub-words) of the expansion.

\begin{lemma}\label{iff}
The $g$-ary expansion of a rational fraction $a/b<1$ in lowest terms has complexity $g^j$ if and only if the following two conditions hold:
\begin{itemize}
\item $0<g^j/b\leq 1$,
\item every interval $[t_s,t_{s+1})$ contains some $r_i\in R$.
\end{itemize}
\end{lemma}
\noindent {\it Proof.} By Lemma \ref{generalizar}, $\{\lfloor tg^j/b\rfloor |t\in B\}=J_g$ if and only if $0<g^j/b\leq 1$, and by Definition \ref{ts}, $\lfloor r_ig^j/b\rfloor=s$ if and only if $r_i\in [t_s,t_{s+1})$.\hfill $\Box$
\begin{definition}\label{gap}
\normalfont Let $r'_0<r'_1<\cdots <r'_{dr -1}$ be the sorted set of power residues $r_i\in R$ generated by $ag^i\equiv r_i\pmod{b}$, and let
\begin{displaymath}
R':=(-1,(r'_0,r'_1,\ldots ,r'_{dr -1}),b).
\end{displaymath}
The {\it gap} of the strictly increasing sequence $R'$ is given by
\begin{displaymath}
G:=\mathrm{max}\{r'_i-r'_{i-1}|i=0,1,\ldots ,dr\},\quad\textrm{where }r'_{-1}:=-1 \textrm{, }r'_{dr}:=b.
\end{displaymath}
\end{definition}
The function gap depends on $g$, $a$ and $b$, where $a/b<1$ is a rational fraction in the lowest terms. For instance, if $g=2$ and $b=43$, we have
\begin{displaymath} 
1/43,  \,\textrm{then } R'=\{-1,1,2,4,8,11,16,21,22,27,32,35,39,41,42,43\}, \,\textrm{thus } G=5;
\end{displaymath}
\begin{displaymath}
3/43,  \,\textrm{then }  R'=\{-1,3,5,6,10,12,19,20,23,24,31,33,37,38,40,43\}, \,\textrm{thus } G=7;
\end{displaymath}
\begin{displaymath}
7/43,  \,\textrm{then } R'=\{-1,7,9,13,14,15,17,18,25,26,28,29,30,34,36,43\},  \,\textrm{thus } G=8.
\end{displaymath}
Since $\varphi (43)=42$ and ord$_2(43)=14$, we have $3$ classes of numerators $a$, i.e., $10/43$ and $3/43$ have the same $R'$. If $g=2$ and $a/b=1/2$, then $G=1$ since $R'=(-1,0,1,2)$.
\begin{theorem}\label{Criterion}
The determination of whether a $g$-ary expansion of a rational fraction $a/b<1$ in the lowest terms has or does not have complexity $g^j$ can be made with the following criterion: 
\begin{displaymath}
\textrm{If }\;\lfloor b/g^j\rfloor \geq G, \textrm{ then complexity }g^j;
\end{displaymath}
\begin{displaymath}
\textrm{If }\;G/2-1<\lfloor b/g^j\rfloor <G, \textrm{ then no decision yet};
\end{displaymath}
\begin{displaymath}
\textrm{If }\;\lfloor b/g^j\rfloor\leq G/2-1,  \textrm{ then no complexity }g^j.
\end{displaymath}
The determination in the second case depends on Lemma \ref{iff}. 
\end{theorem}
\noindent {\it Proof.} We apply Lemma \ref{iff} to the sequence $R'$.

(a) Let $\lfloor b/g^j\rfloor\geq G$. The first condition of Lemma \ref{iff} is fulfilled since $G\geq 1$. Concerning the second condition, we have
\begin{displaymath}
1\leq r'_i-r'_{i-1}\leq G\leq \lfloor b/g^j\rfloor\leq t_{s+1}-t_s\leq\lceil b/g^j\rceil
\end{displaymath}
and therefore,
\begin{displaymath}
r'_i-r'_{i-1}\leq t_{s+1}-t_s\,\; \textrm{for all } i=0,\ldots ,dr \;\; \textrm{and all } s\in J_g.
\end{displaymath}
We choose some fixed $s\in J_g$. The set $\{r'_{i-1}|r'_{i-1}<t_s\}$, with $i=0,\ldots ,dr$, is a non-empty proper subset of $R'$. We define $r'_k:=$max$\{r'_{i-1}|r'_{i-1}<t_s\}$. For the adjacent $r'_{k+1}\in R'$, the following relation holds:
\begin{displaymath}
r'_k<t_s\leq r'_{k+1}<t_{s+1}.
\end{displaymath}
Because, if we assume that $r'_k<t_s<t_{s+1}\leq r'_{k+1}$, then we obtain that $r'_{k+1}-r'_k>t_{s+1}-t_s$ which is impossible. Hence $r'_{k+1}\in [t_s,t_{s+1})$.

(b) We prove that $\lfloor b/g^j\rfloor > G/2-1$ is a necessary condition for $a/b$ to have complexity $g^j$. We assume complexity $g^j$ and therefore, both conditions of Lemma \ref{iff} are fulfilled. 

The master sequence $(t_s)_{s=0}^{g^j}$ has at least three elements since $g\geq 2$ and $j\geq 1$, and we may write
\begin{displaymath}
t_{s-1}\leq r'_{i-1}<t_s\leq r'_i<t_{s+1}.
\end{displaymath}
It follows that $t_{s+1}-t_{s-1}>r'_i-r'_{i-1}\;$ and therefore,
\begin{equation}\label{jose}
t_{s+1}-t_{s-1}=(t_{s+1}-t_s)+(t_s-t_{s-1})>G\quad (s\in J_g).
\end{equation}
Since $1\leq \lfloor b/g^j\rfloor\leq t_{s+1}-t_s\leq \lceil b/g^j\rceil$ for all $s\in J_g$, the left-hand side of the inequality (\ref{jose}) can take only three values, namely: $2\lfloor b/g^j\rfloor$, $2\lfloor b/g^j\rfloor +1$ or $2\lfloor b/g^j\rfloor +2$.
The least bound for $\lfloor b/g^j\rfloor$ is given by $2\lfloor b/g^j\rfloor +2>G$, and hence $\lfloor b/g^j\rfloor > G/2-1.\footnote{We apply Theorem \ref{Criterion} to Stoneham's example, \cite[p. 230]{1969}, $1/59+1/97+1/109$ in base 10. We have $Z=22727$ and $m=623807$. The gap $G$ is $261$. Clearly, we have complexity $g^j$ for $j\leq 3$. For instance, the four-word $9981$ does not occur in the expansion.} $\hfill $\Box$


\section{Frequencies and Transition Matrix}\label{Matrix}
\begin{definition}\label{vec}
\normalfont The {\it frequency vector} for $k\in\mathbb{N}$ is defined by
\begin{displaymath}
V_k:=\big (\nu_k(0),\nu_k(1),\ldots ,\nu_k(s),\ldots ,\nu_k(g^j-1)\big ),
\end{displaymath}
where $\nu_k(s)$ denotes the absolute frequency of the $j$-word $s\in J_g$ in the enlarged repetend $F_g(a/b^k,j)$.
If $\nu_k(s)\not =0$ for all $s\in J_g$, we say that $F_g(a/b^k,j)$ has full complexity $g^j$. Of course, $\sum_{s\in J_g}\nu_k(s)=\textrm{ord}_g(b^k)$ since each of the $di= \textrm{ord}_g(b^k)$ digits in the period is the first digit of some $j$-word in $F_g(a/b^k,j)$.
\end{definition}
\begin{proposition}\label{Dt}
Let $a/b^k$ be a rational fraction in lowest terms, where $\mathrm{gcd}(g,b)=1$ and $k\geq c_g(b)$. Let $j\geq 1$ be a fixed integer and $t\in\mathbb{N}$. Then it holds that
\begin{displaymath}
V_{k+t}=V_k\cdot W_t,
\end{displaymath}
where $W_t=\big(w_{s,s'}\big )$ is a $g^j\times g^j$ matrix --- the transition matrix --- defined by
\begin{displaymath}
w_{s,s'}:= \left\{\begin{array}{ll} \lfloor b^t/g^j\rfloor +1& \textrm{for all $s$ in the sequence}\\ 
&\textrm{$(s'u_t,s'u_t+1,s'u_t+2,\ldots ,s'u_t+u_t-1)\pmod{g^j}$};\\ 
\lfloor b^t/g^j\rfloor& \textrm{for any other $s\in J_g$;} \end{array}\right.
\end{displaymath}
where $u_t$ is defined by $u_t:=b^t \bmod{g^j}$ and $s'\in\{0,1,\ldots , g^j-1\}=J_g$.
\end{proposition}
\noindent {\it Proof.} By Proposition \ref{Xkt}, the list numbers $s'$ of $j$-words in $F_g(a/b^{k+t},j)$ are generated by the expression
\begin{displaymath}
\big (\frac{1}{b^t}(s-n)\big )\bmod{g^j},
\end{displaymath}
where $s$ is a list number of some $j$-word in $F_g(a/b^k,j)$ and $n\in \{0,1,\ldots , b^t-1\}$. We summarize all possibilities in a $g^j\times b^t$ matrix
\begin{displaymath}
E=\big ( e_{s,n}\big )_{s=0,1,\ldots ,g^j-1; n=0,1,\ldots ,b^t-1}
\end{displaymath}
defined by
\begin{displaymath}
e_{s,n}:=\big (\frac{1}{b^t}(s-n)\big )\bmod{g^j}.
\end{displaymath}
The entries of $E$ are all the $s'\in J_g$ that can be generated by Proposition \ref{Xkt} regardless of whether we have $\nu_k(s)=0$ or $\nu_k(s)\not =0$ (see an example in Table \ref{E}).
\begin{table}[htb]
\begin{center}
\begin{tabular}{|c|cccccccccccccc|}
\hline
&n&0&1&2&3&4&5&6&7&8&9&10&11&12\\
\hline
s&&&&&&&&&&&&&&\\
0&&0&2&4&6&8&1&3&5&7&0&2&4&6\\
1&&7&0&2&4&6&8&1&3&5&7&0&2&4\\
2&&5&7&0&2&4&6&8&1&3&5&7&0&2\\
3&&3&5&7&0&2&4&6&8&1&3&5&7&0\\
4&&1&3&5&7&0&2&4&6&8&1&3&5&7\\
5&&8&1&3&5&7&0&2&4&6&8&1&3&5\\
6&&6&8&1&3&5&7&0&2&4&6&8&1&3\\
7&&4&6&8&1&3&5&7&0&2&4&6&8&1\\
8&&2&4&6&8&1&3&5&7&0&2&4&6&8\\
\hline
\end{tabular}
\caption{$s' \equiv \big (\frac{1}{b^t}(s-n)\big )\bmod{g^j}$, where $g=3$, $a/b=7/13$, $j=2$, $t=1$.}
\label{E}
\end{center}
\end{table}
The following properties of $E$ are easy to prove:
\begin{enumerate}
\item $e_{(s+1)\bmod{g^j},\,n}-e_{s,n}\equiv 1/b^t\pmod{g^j}$;
\item $e_{(s+1)\bmod{g^j},\,n+1}\equiv e_{s,n}\pmod{g^j}$\quad for\, $0\leq n<b^t-1$,\\
$e_{(s+1)\bmod{g^j},\,u_t}\equiv e_{s,\,b^t-1}\pmod{g^j}$;
\item $e_{s,\,n+1}-e_{s,n}\equiv -1/b^t\pmod{g^j}$\quad for\, $0\leq n<b^t-1$,\\
$e_{s,\,u_t}-e_{s,\,b^t-1}\equiv -1/b^t\pmod{g^j}$;
\item $e_{(s+u_t)\bmod{g^j},\,n}-e_{s,n}\equiv 1\pmod{g^j}$.
\end{enumerate}

By Property 1, each column is a certain permutation of $(0,1,\ldots ,g^j-1)$. By Property 2, the entries of any diagonal Northwest-Southeast are all the same number. Such a diagonal can be extended by jumping from the bottom of $E$ to row $0$ of the next column, or from the right margin of $E$ to the next row at column $u_t$. By Property 2, Property 3 and because of
\begin{displaymath}
\frac{b^t}{g^j}=\Big\lfloor\frac{b^t}{g^j}\Big\rfloor +\frac{u_t}{g^j},
\end{displaymath}
in each row, a certain permutation of $(0,1,\ldots ,g^j-1)$ is repeated $\lfloor b^t/g^j\rfloor$ times; the remaining $u_t$ numbers appear exactly $(\lfloor b^t/g^j\rfloor +1)$ times in that row. In particular, we have that any $s'\in J_g$ appears exactly $b^t$ times as an entry of $E$.

Let $\mu(s,s')$ denote the number of times a $j$-word $s'\in J_g$ appears in row $s$ of the matrix. This leads to a $g^j\times g^j$ matrix of absolute frequencies $\mu(s,s')$ which will be the claimed {\it transition matrix} $W_t$.

Starting with $e_{0,0}=0$ and applying Property 3 several times, we have the following:

\noindent
since $e_{0,0}=0$,
\begin{displaymath}
\mu(s,0)= \left\{\begin{array}{ll} \lfloor b^t/g^j\rfloor +1& \textrm{for $\;\;0\leq s < u_t$},\\ \lfloor b^t/g^j\rfloor& \textrm{else;} \end{array}\right.
\end{displaymath}
since $e_{u_t,0}=1$,
\begin{displaymath}
\mu(s,1)= \left\{\begin{array}{ll} \lfloor b^t/g^j\rfloor +1& \textrm{for $\;\;u_t\leq s < (2u_t)\bmod{g^j}$},\\ \lfloor b^t/g^j\rfloor& \textrm{else;} \end{array}\right.
\end{displaymath}
finally, since $e_{(s'u_t)\bmod{g^j},0}=s'$ for all $s'\in J_g$,
\begin{displaymath}
\mu(s,s')= \left\{\begin{array}{ll} \lfloor b^t/g^j\rfloor +1& \textrm{for $\;\;(s'u_t)\bmod{g^j}\leq s < \big((s'+1)u_t\big )\bmod{g^j}$},\\ \lfloor b^t/g^j\rfloor& \textrm{else.} \end{array}\right.
\end{displaymath}

Hence the values $\lfloor b^t/g^j\rfloor +1$ form blocks of $u_t$ adjacent numbers in each column $s'$, starting at row $(s'u_t)$ modulo $g^j$ of the matrix $W_t=\mu(s,s')$.

Finally, the multiplication of the vector $V_k$ --- a $1\times g^j$ matrix --- by column $s'$ of $W_t$ yields $\nu_{k+t}(s')$ and therefore, $V_{k+t}=V_k\cdot W_t$ as claimed.\hfill $\Box$
\begin{example}\label{Ex}
\normalfont We calculate the transition matrix $W_t$ of the sequence
\begin{displaymath}
\Big(\frac{7}{13^k}\Big)_{k=1}^{\infty}
\end{displaymath}
in base $g=3$ for $j=2$ and $t=1$. Of course, $a=7$ and $b=13$.

In this case, we have $c_g(b)=c_3(13)=1$.
Indeed, we see that\footnote{We omit the proof. However, there can be surprises. For instance, $\mathrm{ord}_3(26)=3$ and $\mathrm{ord}_3(26^2)=3\cdot 26$; nevertheless, $c_{3}(26)$ is not $1$ but $3$ (see Table \ref{pgm}).}
\begin{displaymath}
\mathrm{ord}_3(13)=3,\quad \mathrm{ord}_3(13^2)=3\cdot 13,\quad \mathrm{ord}_3(13^3)=3\cdot 13^2,\quad \ldots \textrm{ etc.}
\end{displaymath}
Furthermore $u_t=u_1\equiv 13\bmod{3^2}\equiv 4$. Hence the index $s$ goes through the sequence
\begin{displaymath}
(4s', 4s'+1, 4s'+2, 4s'+3) \bmod{3^2}\quad \textrm{for } s'=0, 1,\ldots ,8.
\end{displaymath}
For these $s$ we have $\lfloor b^t/g^j\rfloor +1=\lfloor 13/9\rfloor +1 =2$. So we obtain
\begin{displaymath}
W_1 =
\left( \begin{array}{ccccccccc}
2&1&2&1&2&1&2&1&1\\
2&1&2&1&2&1&1&2&1\\
2&1&2&1&1&2&1&2&1\\
2&1&1&2&1&2&1&2&1\\
1&2&1&2&1&2&1&2&1\\
1&2&1&2&1&2&1&1&2\\
1&2&1&2&1&1&2&1&2\\
1&2&1&1&2&1&2&1&2\\
1&1&2&1&2&1&2&1&2
\end{array}\right)_{.}
\end{displaymath}

The repetend of $7/13$ in base $3$ can be calculated with Lemma \ref{Sg} by choosing $j=1$. We obtain the word $112$. Actually, we have 
\begin{displaymath}
7/13=0.\,112 112 112 112\ldots \,\, .
\end{displaymath}
Hence $F_3(7/13,2)=1121$ is the enlarged repetend. The two-words have the list numbers $4$, $5$ and $7$. Since each of these numbers occurs only once and $J_3=\{0,1,2,3,4,5,6,7,8\}$, we obtain the following frequency vector for $k=1$:
\begin{displaymath}
V_1=(0,0,0,0,1,1,0,1,0).
\end{displaymath}
We apply Proposition \ref{Dt} and get
\begin{displaymath}
V_1\cdot W_1=(0,0,0,0,1,1,0,1,0)\cdot W_1=(3,6,3,5,4,5,4,4,5)=V_2.
\end{displaymath}
Indeed, applying Lemma \ref{Sg} we have
\begin{displaymath}
F_3(7/13^2,2)=0010100120211001021111201222012102122220;
\end{displaymath}
the corresponding list numbers of two-words are the following:
\begin{displaymath}
S_3(7/13^2,2)=(\scriptstyle 0, 1, 3, 1, 3, 0, 1, 5, 6, 2, 7, 4, 3, 0, 1, 3, 2, 7, 4, 4, 4, 5, 6, 1, 5, 8, 8, 6, 1, 5, 7, 3, 2, 7, 5, 8, 8, 8, 6\displaystyle ).
\end{displaymath}
Actually, we have the word $0$, three times; the word $1$, six times; the word $2$, three times; etc. and the word $8$, five times.
\end{example}
\begin{proposition}\label{D1}
Let $W_t$ be as in Proposition \ref{Dt}. Then
\begin{displaymath}
W_t=(W_1)^t.
\end{displaymath}
\end{proposition}
\noindent {\it Proof.} By Proposition \ref{Dt}, we have
\begin{displaymath}
V_{k+1}=V_kW_1,\quad  V_{k+2}=V_{k+1}W_1,\quad  V_{k+3}=V_{k+2}W_1, \ldots \,.
\end{displaymath}
Hence $V_{k+t}=V_k(W_1)^t$. By Proposition \ref{Dt}, also $V_{k+t}=V_kW_t$. Therefore, $W_t=(W_1)^t$.\hfill $\Box$

The sequence
\begin{displaymath}
\big (u_t\big )_{t=1}^\infty =\big (b^t\bmod{g^j}\big )_{t=1}^\infty
\end{displaymath}
is a purely periodic sequence since $b$ and $g$ are relatively prime. Of course, the period length is the order of $b$ modulo $g^j$, that is,
\begin{equation}\label{h0}
t_0:=\textrm{ord}_{b}(g^j).
\end{equation}

The transition matrix $W_t$ can be written as sum
\begin{equation}\label{sum}
W_t=H+W'_t,
\end{equation}
where H has all its entries equal to $\lfloor b^t/g^j\rfloor$, whereas $W'_t$ has only entries $0$ or $1$. Since $b^{t_0}\equiv 1\pmod{g^j}$, we have $u_{t_0}=1$. The entries of the main diagonal Northwest-Southeast of $W_{t_0}$ are all the same number $\lfloor b^{t_0}/g^j\rfloor +1$, and the matrix $W'_{t_0}$ defined by (\ref{sum}) is the identity matrix. Hence $(W'_1, W'_2,\ldots )$ is a periodic sequence with finitely many distinct terms, namely $(W'_1,W'_2,\dots ,W'_{t_0})$.

\begin{proposition}\label{vrel}
Let $a/b^k<1$ be a rational fraction in lowest terms, $\mathrm{gcd}(g,b)=1$ and $k\geq c_g(b)$. Let $j\geq 1$ be a fixed integer and $t \in\mathbb{N}$. The frequency vector $V_{k+t}$ is determined by $V_k$ and $t$ (see an example in Table \ref{tabla}):
\begin{displaymath}
\nu_{k+t}(s')=\Big \lfloor \frac{b^t}{g^j}\Big\rfloor \mathrm{ord}_g(b^k)+\sum_s \nu_k(s),
\end{displaymath}
where $s'\in J_g$ is the list number of some $j$-word in $F_g(a/b^{k+t},j)$, whereas the summation index $s$ goes through the sequence
\begin{displaymath}
(s'u_t, s'u_t+1, s'u_t+2,\ldots , s'u_t+u_t-1) \bmod{g^j}
\end{displaymath}
for $u_t:=b^t \bmod{g^j}$. Clearly,
\begin{displaymath}
0\leq \sum_{s} \nu_k(s)\leq \mathrm{ord}_g(b^t).
\end{displaymath}
\end{proposition}
\noindent {\it Proof.} We have $V_{k+t}=V_k(H+W'_t)=V_kH+V_kW'_t$. The first term of this sum is given by
\begin{displaymath}
V_kH=\big (\lfloor b^t/g^j\rfloor\textrm{ord}_g(b^k),\ldots ,\lfloor b^t/g^j\rfloor\textrm{ord}_g(b^k)\big ),
\end{displaymath}
since $\sum_{s\in J_g} \nu_k(s)=\textrm{ord}_g(b^k)$,
and the second term by
\begin{displaymath}
V_kW'_t=\Big (\sum_{s=0}^{u_t-1} \nu_k(s), \ldots ,\sum_{s=(s'u_t)\bmod{g^j}}^{(s'u_t+u_t-1)\bmod{g^j}} \nu_k(s),\ldots ,\sum_{s=\big((g^j-1)u_t\big)\bmod{g^j}}^{g^j-1}\nu_k(s)\Big ),
\end{displaymath}
since the entries $1$ (respectively $0$) of $W'_t$ correspond to the entries $\lfloor b^t/g^j\rfloor+1$ (respectively $\lfloor b^t/g^j\rfloor$) of $W_t$. Finally, we add the two preceding vectors.\hfill $\Box$

\begin{lemma}\label{vv}
Let $t_0=\mathrm{ord}_b(g^j)$. If $k\geq c_g(b)$, then
\begin{displaymath}
\nu_k(s)-\nu_k(s')=\nu_{k+t_0}(s)-\nu_{k+t_0}(s')
\end{displaymath}
for all $s,s'\in J_g$.
\end{lemma}
\noindent {\it Proof.} Since
\begin{displaymath}
b^{t_0}\equiv 1\pmod{g^j},
\end{displaymath}
we have $u_{t_0}=1$. The entries of the main diagonal Northwest-Southeast of $W_{t_0}$ are all the same number $\lfloor b^{t_0}/g^j\rfloor +1$, and the matrix $W'_{t_0}$ of Equation (\ref{sum}) is the identity matrix. By Proposition \ref{vrel}, we have
\begin{displaymath}
\nu_{k+t_0}(s)=\Big\lfloor\frac{b^{t_0}}{g^j}\Big\rfloor \textrm{ord}_g(b^k)+\nu_k(s),
\end{displaymath}
and the same holds for $s'$ instead of $s$.\hfill $\Box$

The following proposition shows that the oscillation of absolute frequencies of $j$-words in $F_g(a/b^k,j)$ is bounded for all $k\in\mathbb{N}$.
\begin{proposition}\label{cota}
Let $a/b$ be a purely periodic proper fraction, $t_0=\mathrm{ord}_b(g^j)$ and $j,k\in \mathbb{N}$. There is a nonnegative integer $C_g(a/b,j)$ such that
\begin{displaymath}
\mid \nu_k(s)-\nu_k(s')\mid \,\leq C_g(a/b,j)
\end{displaymath}
for all $s,s'\in J_g$ and all $k\geq c_g(b)$, where $c_g(b)$ is the integer of Definition \ref{cg}.
\end{proposition}
\noindent {\it Proof.} Indeed, by Lemma \ref{vv},
\begin{displaymath}
C_g(a/b,j):=\textrm{max}\big\{|\nu_{c_g(b)+t}(s)-\nu_{c_g(b)+t}(s')|\,\big |\, s,s'\in J_g,\,\, t\in\{0,1,\ldots ,t_0-1\}\big\}
\end{displaymath}
is the claimed $C_g(a/b,j)$.\hfill $\Box$

Note that $C_g(a/b,j)$ is a bound for all $k\geq c_g(b)$. The oscillation corresponding to the finitely many $k<c_g(b)$ might exceed this bound.

\begin{example}\label{Ex2}
\normalfont We compute $C_3(7/13,2)$ and, as shown in example 1, the transition matrix is $W_1$. In Table \ref {tabla} we have the frequency vectors for $k\in\{1,2,3,4,5\}$.
\begin{table}[htb]
\begin{center}
\begin{tabular}{|r@{}ccccccccc|}
\hline
$V_1=$&(0,&0,&0,&0,&1,&1,&0,&1,&0);\\
$V_2=$&(3,&6,&3,&5,&4,&5,&4,&4,&5);\\
$V_3=$&(56,&56,&56,&57,&57,&56,&55,&57,&57);\\
$V_4=$&(732,&732,&732,&732,&733,&733,&732,&733,&732);\\
$V_5=$&(9519,&9522,&9519,&9521,&9520,&9521,&9520,&9520,&9521).\\
\hline
\end{tabular}
\caption{Frequency vectors of $7/13^k$ in base $g=3$ for $j=2$ and $k\in\{1,2,3,4,5\}$.}
\label{tabla}
\end{center}
\end{table}

\noindent For instance, in the enlarged repetend $F_3(7/13^3,2)$ (third row), the two-word $11$ (having list number $4$) occurs $57$ times, and the two-word $12$ occurs $56$ times, meaning that we will have the following oscillations of absolute frequencies: $\sigma_3(7/13,2)=1$, $\sigma_3(7/13^2,2)=3$ and $\sigma_3(7/13^3,2)=2$, with all other oscillations for $k>3$ being a repetition of these values. Hence $C_3(7/13,2)=3$. 
\end{example}

\begin{proposition}\label{antespure}
Let $a/b^{k+t}<1$ be a rational fraction in lowest terms, $\mathrm{gcd}(g,b)=1$ and $t\in\mathbb{N}$. If $k\geq c_g(b)$,  then
\begin{displaymath}
\Big |\nu_{k+t}(s)- \frac{b^t}{g^j}\mathrm{ord}_g(b^k)\Big |< \mathrm{ord}_g(b^k)
\end{displaymath}
for any $j$-word $s$ in the enlarged period $F_g(a/b^{k+t},j)$.
\end{proposition}
\noindent {\it Proof.} Let $\mu:=\mathrm{ord}_g(b^k)$. By Proposition \ref{vrel}, we have
\begin{equation*}
b^t\mu/g^j-\mu < \lfloor b^t/g^j\rfloor \mu +0\leq \nu_{k+t}(s)\leq \lfloor b^t/g^j\rfloor\mu +\mu < b^t\mu/g^j +\mu.\\[-20pt]
\end{equation*}\hfill $\Box$

\begin{theorem}\label{pure}
For any purely periodic proper fraction $a/b^{k+t}$ with $k\geq c_g(b)$, $t\in\mathbb{N}$, for all $j\geq 1$ each sequence of $j$ digits occurs in the $g$-ary repetend of $a/b^{k+t}$ with a relative frequency that approaches $1/g^j$ for an increasing $t$.
\end{theorem}
\noindent {\it Proof.} By Proposition \ref{antespure} and Definition \ref{cg} (definition of $c_g(b)$\,), the statement holds for the enlarged repetend of $a/b^{k+t}$. Indeed,
\begin{displaymath}
\Big | \frac{\nu_{k+t}(s)}{\mathrm{ord}_g(b^{k+t})} - \frac{1}{g^j}\Big |\leq \frac{1}{b^t}\,.
\end{displaymath}

When leaving out the last $j-1$ digits of $F_g(a/b^{k+t},j)$, the new absolute frequencies will be
\begin{displaymath}
\nu'_{k+t}(s)=\nu_{k+t}(s)-i \quad \textrm{for some }\, i\in \{0,1,\ldots , j-1\}.
\end{displaymath}
Furthermore,
\begin{displaymath}
\sum_{s\in J_g} \nu'_{k+t}(s)=\textrm{ord}_g(b^{k+t})-(j-1).
\end{displaymath}

The new relative frequency
\begin{displaymath}
\nu'_{\textrm{rel}}(s)=\frac{\nu_{k+t}(s)-i}{\textrm{ord}_g(b^{k+t})-(j-1)}
\end{displaymath}
approaches $\nu_{k+t}(s)/\textrm{ord}_g(b^{k+t})$ for increasing $t$.\footnote{ Complexity $g^j$ and ($j,\varepsilon $)-normality, Stoneham \cite{1969}, of an individual rational fraction mean the same thing. Indeed, if $Z/b$ is ($j,\varepsilon $)-normal in base $g$, then $j\geq 1$, and the word $B_j$ occurs at least once. Conversely, if we have complexity $g^j$ for some fixed $j=j_0\geq 1$, then every $B_j$ occurs at least once, and we have ($j,\varepsilon $)-normality for $j=j_0$ and for all $\varepsilon $ larger than each of the finitely many $|N(B_j,g)/\omega (b)-1/g^j|$. }\hfill $\Box$


\section{Generalization for Proper Fractions}\label{Prop}
A $g$-ary expansion of a rational number
\begin{displaymath}
a/b=x_{-\ell}\cdots x_{-1}x_0.\,x_1x_2\cdots \quad (\ell\in \mathbb{N}\cup  \{0\})
\end{displaymath}
is eventually periodic. There is a positive integer $n\in\mathbb{N}$ and an index $i_0\in\mathbb{N}$ such that
\begin{displaymath}
x_i=x_{i+n} \quad \textrm{for all }i\geq i_0.
\end{displaymath}
If $n$ and $i_0$ are the smallest positive integers with this property, we say that
\begin{itemize}
\item $n\in \mathbb{N}$ is the {\it period};
\item the word $x_1x_2\cdots x_{i_0-1}$ is the {\it transient};\footnote{For $i_0=1$, the transient is the empty word.}
\item the word $x_{i_0}x_{i_0+1}\cdots x_{i_0+n-1}$ is the {\it repetend};
\item the word $(x_{i_0}x_{i_0+1}\cdots x_{i_0+n-1}) (x_{i_0}x_{i_0+1}\cdots x_{i_0+j-2})$ is the {\it enlarged repetend}, where $j\geq 1$ is an integer.\footnote{For $j=1$, the second factor is the empty word.}
\end{itemize}

P. Bundschuh  \cite{1988} has shown that in the $g$-adic expansion of a rational number $a/b$ in the lowest terms, the length of the transient and the length of the repetend depend only on $b$ and $g$, but not on $a$. This is what our following lemma proves.

\begin{lemma}\label{trans}
Let $a/b$ be a proper fraction in lowest terms and $\{p_1,p_2,\ldots , p_i, \ldots ,p_t\}$ the non-empty set of primes dividing simultaneously the base $g\geq 2$ and the denominator $b\in\mathbb{N}\setminus \{1\}$. Define
\begin{equation*}
g=g'\prod_{i=1}^t p_i^{g_i}\quad\textrm{and}\quad b=b'\prod_{i=1}^t p_i^{b_i},
\end{equation*} 
where $\mathrm{gcd}(g,b')=1$, and where $g_i$ and $b_i$ are the highest exponents of $p_i$ in the prime factorizations of $g$ and $b$.
The length $f$ of the transient in $a/b=0.\,x_1x_2\cdots$ is given by
\begin{equation*}
f=\mathrm{min}\big\{m\in\mathbb{N}\,\mid \, mg_i\geq b_i \quad\textrm{for all }i\in\{1,2,\ldots ,t\}\big\}.\footnote{P. Bundschuh \cite{1988} defined $b^*:=b'$, $b^{**}:=\prod_{i=1}^t p_i^{b_i}$, ord$_{b^*}g:=$ ord$_g(b')$, $\mathrm{Min}\big\{\mu\in\mathbb{N}_0\,:\,b^{**}\mid \, g^\mu \big\}:=f$.}
\end{equation*} 
The repetend of $a/b$ is that of the purely periodic proper fraction
\begin{equation*}
\frac{\Big (a (g')^f\prod_{i=1}^t p_i^{\,fg_i-b_i}\Big )\bmod{b'}}{b'}
\end{equation*} 
provided that $b'\not =1$.
\end{lemma}
\noindent {\it Proof.} Indeed, we have
\begin{equation}\label{calc}
\frac{a}{b}=\frac{a}{b'\prod p_i^{b_i}}\cdot \frac{(g')^f\prod p_i^{\,fg_i-b_i}}{(g')^f\prod p_i^{\,fg_i-b_i}}=\frac{1}{g^f}\cdot \frac{a(g')^f \prod p_i^{\,fg_i-b_i}}{b'}
\end{equation}
and therefore,
\begin{displaymath}
\frac{a(g')^f \prod p_i^{\,fg_i-b_i}}{b'}=x_1x_2\cdots x_f.\,x_{f+1}x_{f+2}\cdots ,
\end{displaymath}
where the leading digits $x_1, x_2,\ldots$ may be zeros. Because of gcd$(g,b')=1$ and $b'\not =1$, the expansion $0.\,x_{f+1}x_{f+2}\cdots$ is purely periodic with period ord$_g(b')$, and the word $x_1x_2\cdots x_f$ is the transient of $a/b$.\hfill $\Box$

If $b'=1$, then all prime factors of $b$ divide $g$. In this case, the $g$-ary expansion of $a/b$ is eventually constant and equal to $0$ or $g-1$.
\begin{lemma}\label{kf}
Let $a/b$ be a proper fraction in the lowest terms, $k\in\mathbb{N}$ and $b'\not =1$ as in Lemma \ref{trans}. The repetend of $a/b^k$ is that of the purely periodic proper fraction
\begin{equation*}
\frac{\Big (a(g')^{kf} \prod_{i=1}^t p_i^{\,k(fg_i-b_i)}\Big )\bmod{(b')^k}}{(b')^k}.
\end{equation*} 
\end{lemma}
\noindent {\it Proof.} Apply Equation (\ref{calc}).\hfill $\Box$

Lemma \ref{kf} shows that the repetends of proper fractions $a/b^k$ are the repetends of purely periodic proper fractions 
of the form
\begin{displaymath}
\frac{a'_k}{(b')^k},
\end{displaymath}
where the numerators $a'_k$ depend on $k$, and where $b'\not =1$ is the $g$-free part of the denominator $b$.

\begin{proposition}\label{seq}
Let $b$ be a positive integer with at least one prime factor that does not divide the base $g\geq 2$.
Let $\big(a_k\big)_{k=1}^\infty$ be a sequence of positive integers relatively prime to $b$.

For all $j\geq 1$, each sequence of $j$ digits occurs in the enlarged repetend of $a_k/b^k$ with a relative frequency that approaches $1/g^j$ for an increasing $k$.
\end{proposition}
\noindent {\it Proof.} The difference $a_k/b^k-\big(a_k\bmod{b^k}\big)/b^k$ is a nonnegative integer. Hence $a_k/b^k$ and $\big(a_k\bmod{b^k}\big)/b^k$ have the same transient and also the same repetend. Without loss of generality, we assume that $a_k/b^k$ are already proper fractions and consider the sequence 
\begin{displaymath}
P:=\Big(\frac{a_k}{b^k}\Big)_{k=1}^\infty.
\end{displaymath}

Of course, the $g$-free part $b'$ of $b$ is not $1$. We apply Lemma \ref{kf} on each term of $P$. We get a sequence
\begin{displaymath}
P':=\Big(\frac{a'_k}{(b')^k}\Big)_{k=1}^{\infty}
\end{displaymath}
of purely periodic proper fractions. The corresponding terms of $P$ and $P'$ have the same repetends.

We now consider the new sequence
\begin{displaymath}
P'':=\bigg(\frac{a'_k \bmod{\big((b')^{c_g(b')} \big)}}{(b')^k}\bigg)_{k=1}^\infty:=\bigg(\frac{a''_k}{(b')^k}\bigg)_{k=1}^{\infty}.
\end{displaymath}
The first $dn=c_g{(b')}$ numerators of $P'$ remain unchanged. Let $k=c_g(b')+t$ ($t\in\mathbb{N})$. We claim that sequences
\begin{displaymath}
R_g\Big(\frac{a'_k}{(b')^k},j\Big)=\big(r'_v\big)_{v=0}^{\textrm{ord}_g((b')^k)-1} \quad\textrm{and}\quad R_g\Big(\frac{a''_k}{(b')^k},j\Big)=\big(r''_u\big)_{u=0}^{\textrm{ord}_g((b')^k)-1},
\end{displaymath}
except for left rotations, have the same remainders.

Indeed, the two sequences share the remainder $a''_k$. Of course, $a''_k=r''_0$. Furthermore we have
\begin{displaymath}
R_g\Big(\frac{a'_{c_g(b')+t}}{(b')^{c_g(b')}},j\Big)=R_g\Big(\frac{a''_{c_g(b')+t}}{(b')^{c_g(b')}},j\Big):=\big(r_i\big).
\end{displaymath}
Note that this equality holds even if $a'_{c_g(b')+t}>(b')^{c_g(b')}$. We apply Lemma \ref{kt} and get the representation 
\begin{displaymath}
a'_{c_g(b')+t}=r'_0=r_{i_0}+n_0(b')^{c_g(b')}
\end{displaymath}
for an appropriate $r_{i_0}\in \big(r_i\big)$ and some $n_0\in\{0,1,\ldots ,(b')^t-1\}$.
Since $a''_{c_g(b')+t}\equiv a'_{c_g(b')+t}\pmod{(b')^{c_g(b')}}$, we have $r_{i_0}= a''_{c_g(b')+t}$. Therefore,
\begin{displaymath}
r'_0=a''_{c_g(b')+t}+n_0(b')^{c_g(b')}.
\end{displaymath}
By Lemma \ref{kt}, there is a remainder $r'_{v_0}\in\big(r'_v\big)$ such that
\begin{displaymath}
r'_{v_0}=a''_{c_g(b')+t}+ 0\cdot (b')^{c_g(b')}=a''_{c_g(b')+t}.
\end{displaymath}
Hence $a''_k$ also belongs to the sequence $\big(r'_v\big)$.

Since $\big(r'_v\big)$ and $\big(r''_u\big)$ are left rotations of the same sequence, the corresponding sequences of $j$-words also are left rotations of the same sequence by Lemma \ref{rs}. It follows that the corresponding terms of $P'$ and $P''$, except for left rotations, have the same repetends. 

Let $D_{c_g(b')}$ be the cyclic subgroup of $(\mathbb{Z}/(b')^{c_g(b')}\mathbb{Z})^\times $ generated by $g$ modulo $(b')^{c_g(b')}$. Furthermore, let $\{z_1,\ldots , z_h,\ldots ,z_{\ell}\}$ be representatives of the
\begin{displaymath}
\ell=\frac{\varphi\big((b')^{\textrm{ord}_g(b')}\big)}{\textrm{ord}_g((b')^{\textrm{ord}_g(b')}\big)}
\end{displaymath}
cosets of $D_{c_g(b')}$.

Starting at $k=c_g(b')$, we replace each numerator $a''_k$ of $P''$ by the representative $z_h$ of the corresponding coset $z_hD_{c_g(b')}$ and get the sequence
\begin{displaymath}
P''':=\Big(\frac{y_k}{(b')^k}\Big)_{k=c_g(b')}^\infty,
\end{displaymath}
where $y_k\in\{z_1,\ldots, z_{\ell}\}$. The corresponding terms of $P$ and $P'''$ have left rotations of the same repetends. 

The sequence $P'''$ splits into finitely many subsequences $\big(z_h/(b')^{k_i}\big)_{k_i}$ of terms with equal numerator $z_h$. Any such subsequence is also a subsequence of
\begin{displaymath}
Z_h:=\Big(\frac{z_h}{(b')^k}\Big)_{k=c_g(b')}^\infty,
\end{displaymath}
and the statements of Section \ref{Matrix} hold for this sequence.

In particular, the oscillation of absolute frequencies of $j$-words in $F_g(z_h/(b')^k,j)$ is bounded by $C_g(z_h/b',j)$ for all $k\geq c_g(b')$ (see Proposition \ref{cota}). Since there are only finitely many sequences $Z_h$,
\begin{displaymath}
c:=\textrm{max}\big\{C_g(z_h/b',j)\,\big |\, z_h\in\{z_1,\ldots, z_{\ell}\}\big\}
\end{displaymath}
is a bound for all $Z_h$. Hence the oscillation of absolute frequencies of $j$-words in the enlarged repetend of the proper fractions $a_k/b^k$ is bounded too for all $k\in\mathbb{N}$.\hfill $\Box$
\begin{theorem}\label{proper}
Let $b>1$ be a positive integer and $\big(a_k\big)_{k=1}^\infty$ a sequence of positive integers relatively prime to $b$. For all $j\geq 1$, each sequence of $j$ digits occurs in the $g$-ary repetend of $a_k/b^k$ with a relative frequency that approaches $1/g^j$ for an increasing $k$, unless all prime factors of $b$ divide the base $g\geq 2$.
\end{theorem}
\noindent {\it Proof.} The statement holds for the enlarged repetend by Proposition \ref{seq}. It also holds for the repetend by the same arguments as in the proof of Theorem \ref{pure}. 

If all prime factors of $b$ divide $g$, then the $g$-ary expansion of $a/b^k$ is eventually constant and equal to $0$ or $g-1$ for all $k\in\mathbb{N}$.\hfill $\Box$


\section{A Helpful Sequence}
\begin{lemma}\label{viejo}
Let $a/b\in \mathbb{Q}$ be a positive fraction in lowest terms: $a\in \mathbb{N}$, $b \in \mathbb{N}\setminus\{1\}$, and $\emph{gcd}(a,b)=1$.
Let $A:=\{0,1,2,\ldots ,a-1\}$ and $B:=\{0,1,2,\ldots ,b-1\}$. Then the following statements hold.
\begin{enumerate}
\item $\Big\{\big\lfloor t\cdot \frac{a}{b}\big\rfloor \,\Big |\, t\in B\Big\} = A$ \quad if and only if $\quad 0<\frac{a}{b}<1$.
\item Let $0<\frac{a}{b}<1$ and $u:=b \bmod{a}$. In the sequence
\begin{displaymath}
\bigg (\Big \lfloor t\cdot\frac{a}{b} \Big \rfloor \bigg )_{t=0}^{b-1},\\
\end{displaymath}
each element of $A$ occurs either $\lceil b/a \rceil$ or $\lfloor b/a \rfloor$ times. There are exactly $u$ elements of $A$ occurring $\lceil b/a \rceil$ times; each of the remaining $a-u$ elements occurs $\lfloor b/a \rfloor$ times.
\end{enumerate}
\end{lemma}
\noindent {\it Proof}. Lemma \ref{viejo} is a simplified version of Proposition 23 in \cite{2012}.\hfill $\Box$
\begin{lemma}\label{generalizar}
Let $z/m>0$ be a rational fraction such that $\mathrm{gcd}(z,m)=d$.
Let $Z=\{0,1,\ldots ,z-1\}$ and $M=\{0,1,\dots ,m-1\}$. Then the following statements hold.
\begin{enumerate}
\item $\Big\{\big\lfloor t\cdot \frac{z}{m}\big\rfloor\big | t\in M\Big\}=Z$ \quad if and only if \quad $0<\frac{z}{m}\leq 1$;
\item Let $0<\frac{z}{m}<1$ be and $u:=\frac{m}{d}\bmod{\frac{z}{d}}$. In the sequence
\begin{displaymath}
S:=\bigg( \Big\lfloor t\cdot \frac{z}{m}\Big\rfloor \bigg )_{t=0}^{m-1},
\end{displaymath}
each element of $Z$ occurs either $\lceil m/z\rceil$ or $\lfloor m/z\rfloor$ times. There are exactly $ud$ elements of $Z$ occurring $\lceil m/z\rceil$ times, and $0$ is one of them; each of the remaining $z-ud$ elements occurs $\lfloor m/z\rfloor$ times.
\end{enumerate}
\end{lemma}
\noindent {\it Proof.} Since $z=da$, $m=db$ and gcd$(a,b)=1$, we have
\begin{displaymath}
Z=\bigcup_{x=0}^{d-1}\{xa,xa+1,\ldots ,xa+a-1\},\quad
M=\bigcup_{x=0}^{d-1}\{xb,xb+1,\ldots ,xb+b-1\}.
\end{displaymath}
Each $t\in M$ has a unique representation
\begin{displaymath}
t=xb+y, \quad\textrm{where}\quad 0\leq y\leq b-1 \quad\textrm{and}\quad 0\leq x\leq d-1.
\end{displaymath}
It follows that
\begin{displaymath}
\Big\lfloor t\cdot\frac{z}{m}\Big\rfloor = xa+\Big\lfloor y\cdot \frac{a}{b}\Big\rfloor.
\end{displaymath}
Hence
\begin{displaymath}
S=(S',a+S',2a+S',\ldots ,(d-1)a+S'),\quad\textrm{where}\quad S':=\bigg(\Big\lfloor y\cdot\frac{a}{b}\Big\rfloor\bigg)_{y=0}^{b-1}.
\end{displaymath}
Of course, $a+S'$ means that we add $a$ to each term of $S'$.

We apply Lemma \ref{viejo} to the sequence $S'$. The set of integers in $S'$ is $A:=\{0,1,\ldots ,a-1\}$ if and only if $0<a/b<1$, and $u$ elements of $A$ occur $\lceil b/a\rceil$ times and  $a-u$ elements $\lfloor b/a\rfloor$ times in $S'$. The number $0$ occurs $\lceil b/a\rceil$ times, because it is the only integer in the sequence $(ya/b)_{y=0}^{b-1}$. If $A$ is the set of integers in $S'$, with $0<a/b<1$, then $S$ is a sequence of consecutive integers from $0$ to $da-1$, wherein $ud$ integers are repeated $\lceil m/z\rceil$ times and $z-ud$ elements $\lfloor m/z\rfloor$ times. Therefore, $Z$ is the set of integers in $S$, and this is also true for $z/m=1$ since then $Z=M$.\hfill $\Box$

\end{document}